\documentclass{amsart}%article

\usepackage{amssymb,amsfonts,mathrsfs, pstricks,epsf}
\usepackage[dvips]{graphicx}

\theoremstyle{plain}
\newtheorem{thm}{Theorem}[section]
\newtheorem{prop}[thm]{Proposition}
\newtheorem{lem}[thm]{Lemma}

\newtheorem{conj}[thm]{Conjecture}
\newtheorem{centralconj}[thm]{Central~Conjecture}

\theoremstyle{definition}

\theoremstyle{remark}
\newtheorem{rem}[thm]{Remark}

\numberwithin{equation}{section}

\newcommand{\Rth}{\textbf{R}$^3$}

\newcommand{\Rf}{\textbf{R}$^4$}

\newcommand{\Rn}{\textbf{R}$^n$}
\newcommand{\RN}{\textbf{R}$^N$}
\newcommand{\mRN}{\mathbf{R}^N}

\newcommand{\Ton}{\textbf{T}$^1$}
\newcommand{\Ttw}{\textbf{T}$^2$}
\newcommand{\Tth}{\textbf{T}$^3$}
\newcommand{\mTth}{\mathbf{T}^3}

\newcommand{\sixty}{\frac{\pi}{3}}

\newcommand{\onetwenty}{120^{\circ}}

\newcommand{\vol}{\mathrm{vol}}
\newcommand{\area}{\mathrm{area}}

\newcommand{\SE}{\emph{Surface Evolver}}

\begin{document}
\title{Double Bubbles in the 3-Torus}
\date{\today}

\author[M. Carri\'on~\'Alvarez]{Miguel Carri\'on~\'Alvarez}
\address{University of California -- Riverside}
\email{miguel@math.ucr.edu}

\author[J. Corneli]{Joseph Corneli}
\address{New College of Florida}
\email{jcorneli@virtu.sar.usf.edu}

\author[G. Walsh]{Genevieve Walsh}
\address{University of California -- Davis}
\email{gwalsh@math.ucdavis.edu}

\author[S. Beheshti]{Shabnam Beheshti}
\address{Texas Tech University }
\email{amper@hilbert.math.ttu.edu}

\thanks{Prepared on behalf of the participants in the Clay Mathematics
Institute Summer School on the Global Theory of Minimal Surfaces, held
at the Mathematical Sciences Research Institute in Berkeley,
California, Summer 2001.}

\begin{abstract}
We present a conjecture, based on computational results, on the area
minimizing way to enclose and separate two arbitrary volumes in the
flat cubic three-torus \Tth. For comparable small volumes, we prove
that an area minimizing double bubble in \Tth\ is the standard double
bubble from \Rth.
\end{abstract} 

\maketitle

\section{Introduction}
Our Central Conjecture 2.1 states that the ten different types of
two-volume enclosures pictured in Figure 1 comprise the complete set
of surface area minimizing \emph{double bubbles} in the flat cubic
three-torus \Tth. Our numerical results, summarized in Figure
\ref{phase_portrait}, indicate the volumes for which we conjecture
that each type of double bubble minimizes surface area. Our main
theorem, Theorem 4.1, states that given any fixed ratio of volumes,
for small volumes, the minimizer is the standard double bubble. This
result applies to any smooth flat Riemannian manifold of dimension
three or four with compact quotient by its isometry group.

%*CATALOG
\begin{figure}[h]
\begin{center}
\begin{tabular}{@{\extracolsep{1 cm}} cc}
        \multicolumn{2}{c}{\ }\\  
        \multicolumn{2}{c}{\ }\\
        \multicolumn{2}{c}{\ }\\

    %SDB
    \rput{0}(0.2,1.9){
    \epsfysize=.75in
        \epsffile{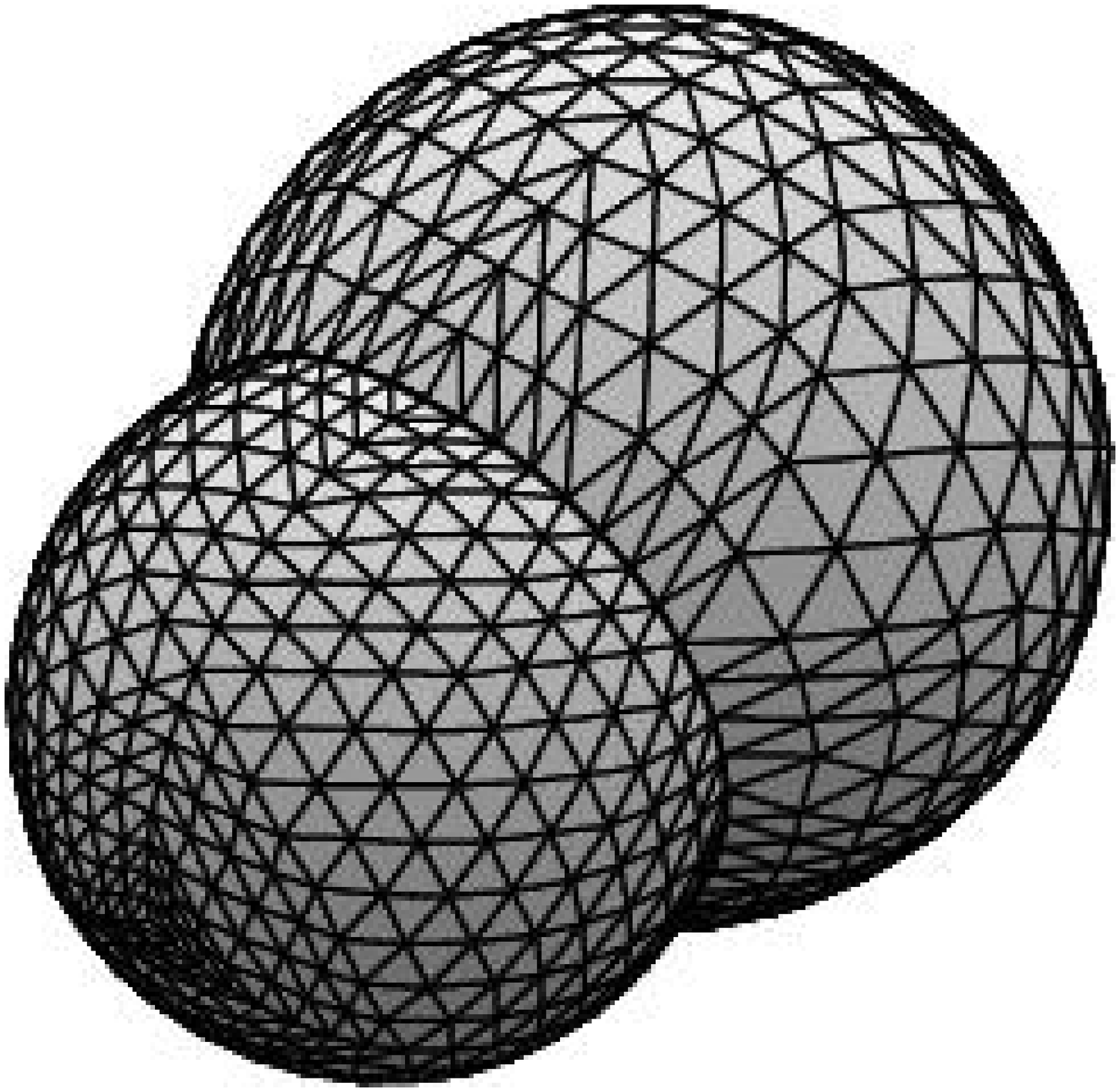}
        }
    \pspicture(-1,-.75)(1,.75)
    \rput*{0}(-1,-.5){ \begin{scshape} Standard Double Bubble \end{scshape}}
        \endpspicture

        &

        %Delauney Chain
    \rput{-30}(2.,1.9){
    \epsfysize=1.0in
        \epsffile{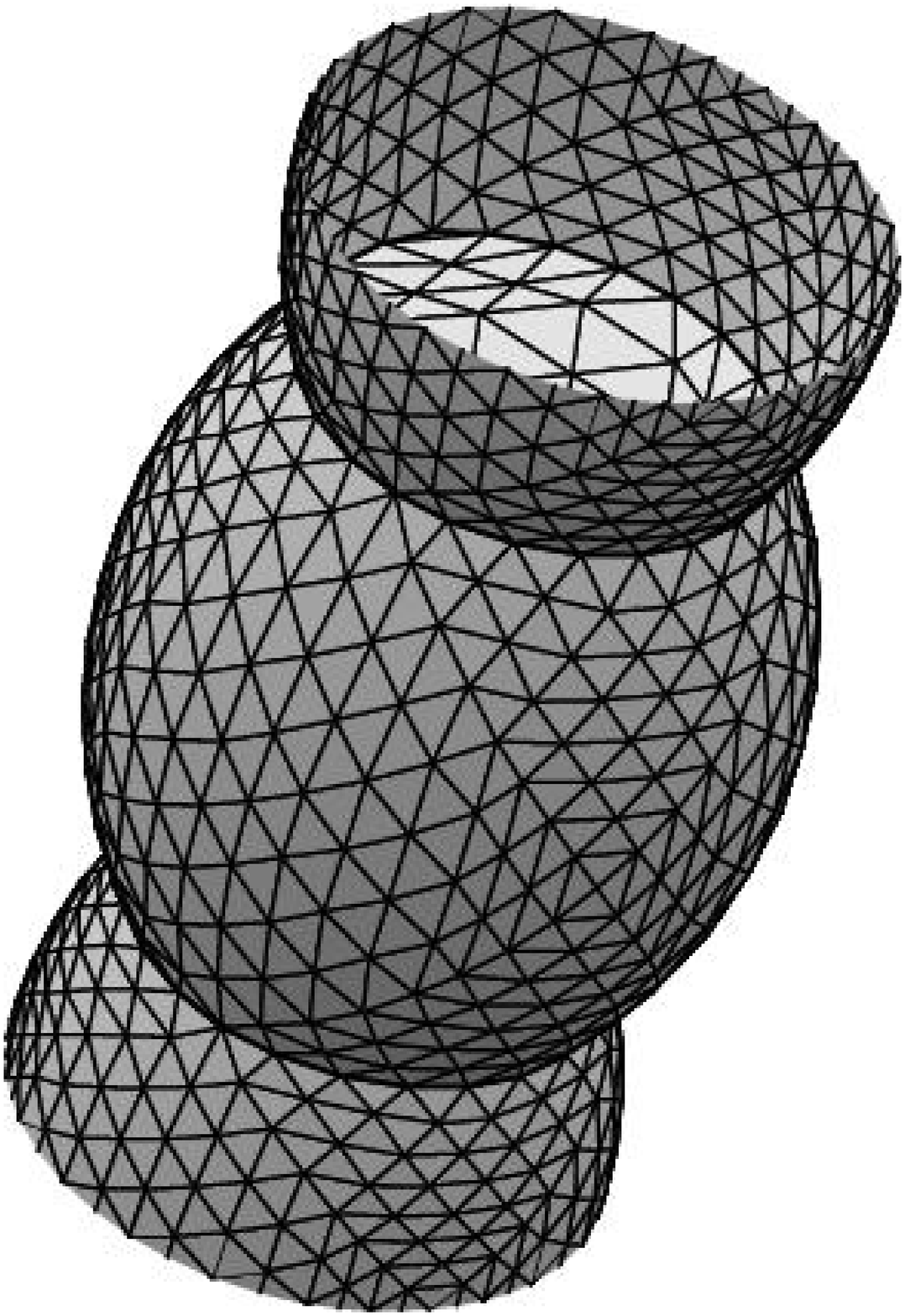}
        }
        \pspicture(-1,-.75)(1,.75)
        \rput*{0}(.8,-.5){ \begin{scshape} D\'elauney Chain \end{scshape}}
        \endpspicture

    \\ 

        \multicolumn{2}{c}{\ }\\  
        \multicolumn{2}{c}{\ }\\  
        \multicolumn{2}{c}{\ }\\
        \multicolumn{2}{c}{\ }\\
        \multicolumn{2}{c}{\ }\\

\end{tabular}

\begin{tabular}{@{\extracolsep{1 cm}} ccc}

    %Cylinder Lens
    \rput{0}(0.,1.9){
    \epsfysize=.80in
        \epsffile{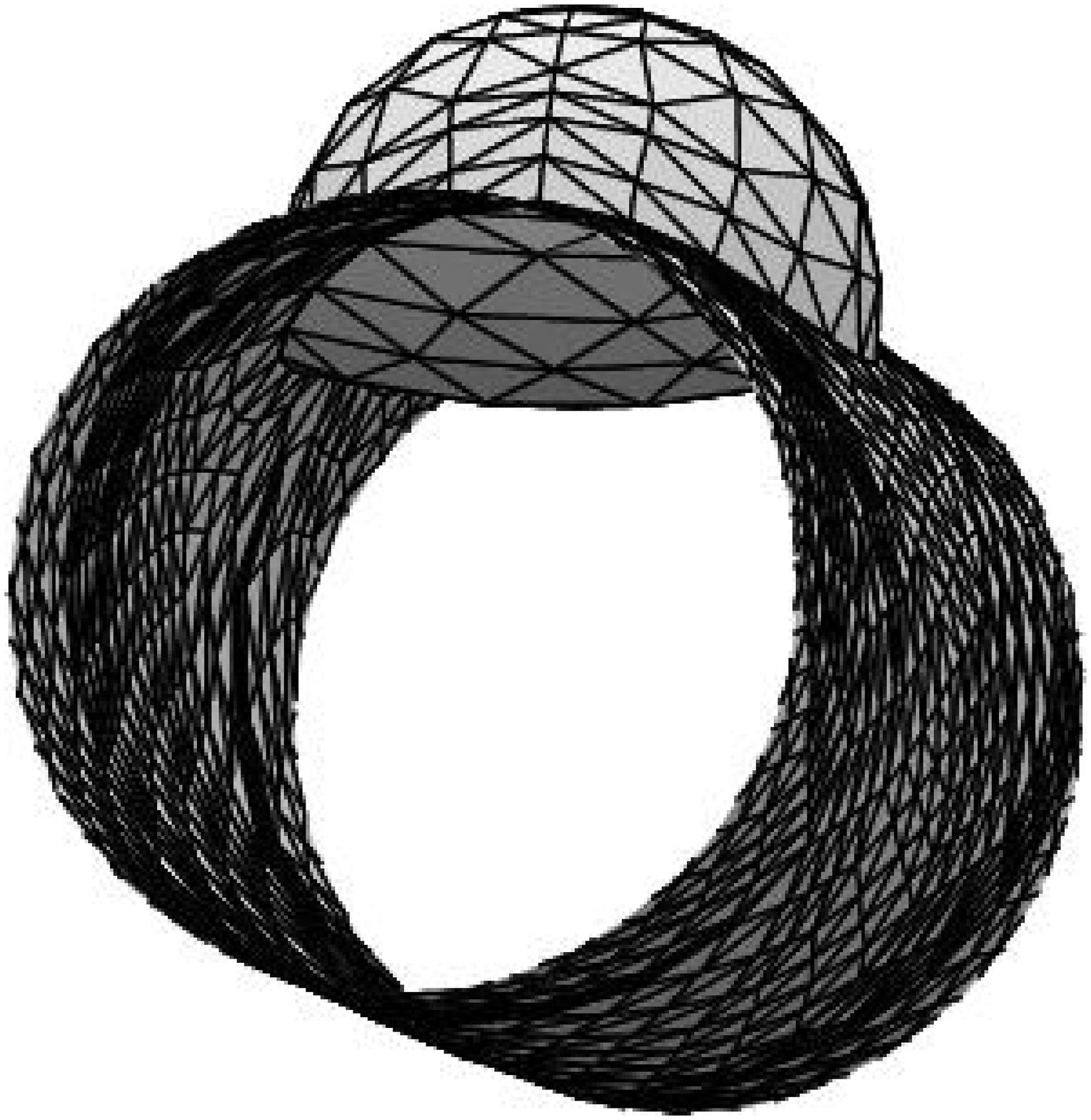}
        }
    \pspicture(-1,-.75)(1,.75)
    \rput*{0}(-1.,-.5){ \begin{scshape} Cylinder Lens \end{scshape}}
        \endpspicture

        &

        %Cylinder Cross
    \rput{0}(1.,1.9){
    \epsfysize=.80in
        \epsffile{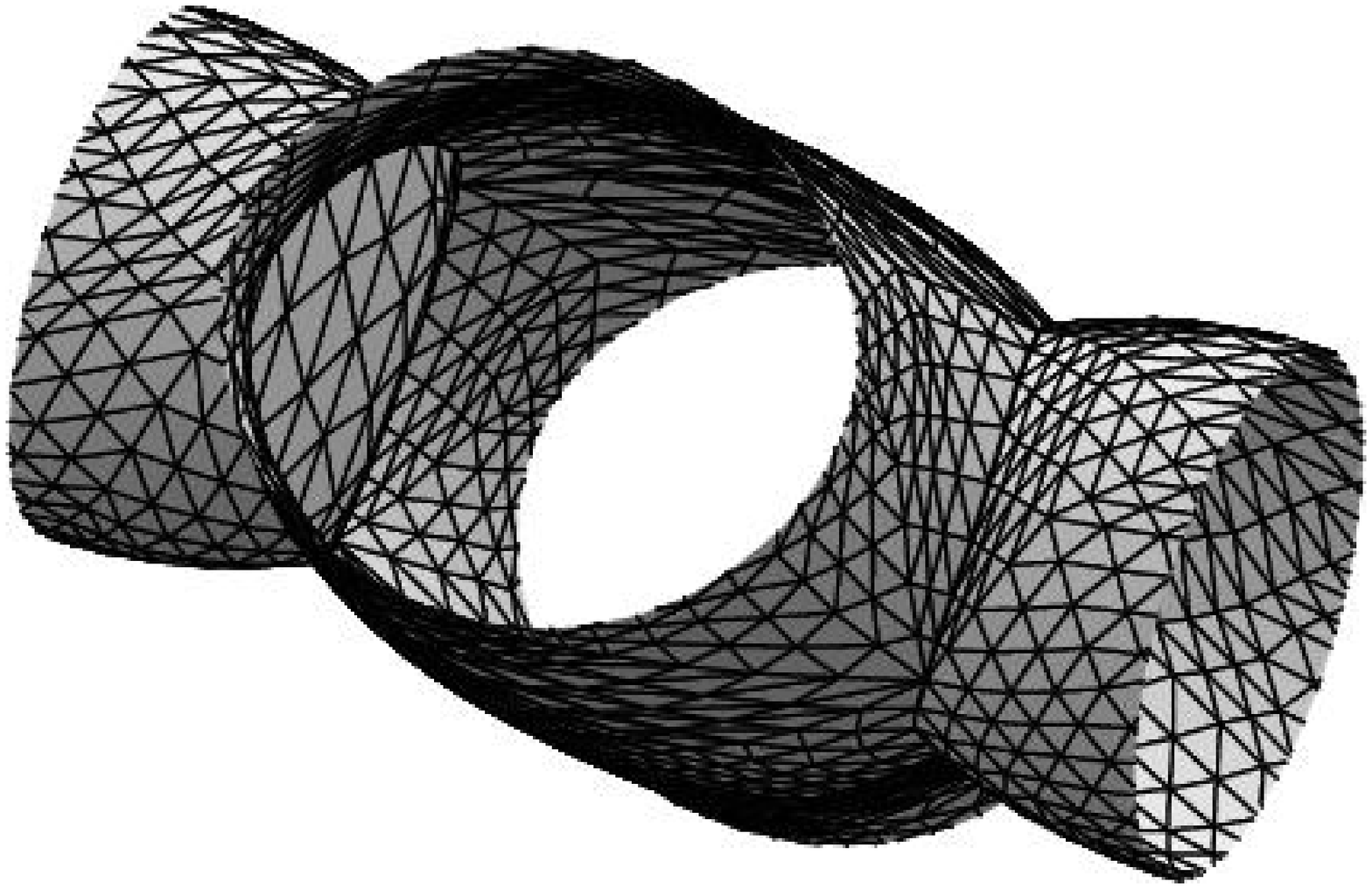}
        }
        \pspicture(-1,-.75)(1,.75)
    \rput*{0}(0.,-.5){ \begin{scshape} Cylinder Cross \end{scshape}}
        \endpspicture

        &

        %Double Cylinder
    \rput{0}(1.6,1.9){
    \epsfysize=.80in
        \epsffile{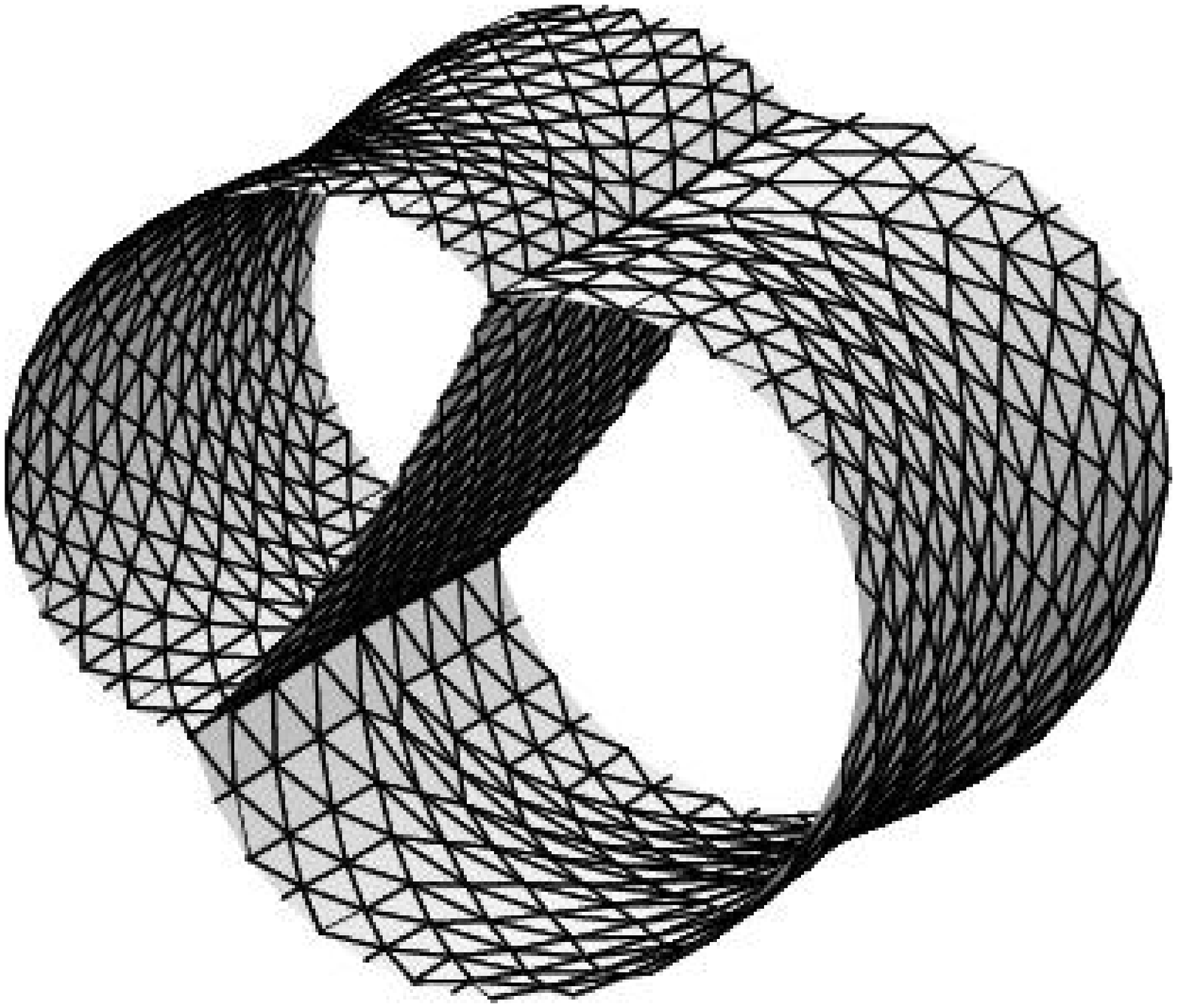}
        }
        \pspicture(-1,-.75)(1,.75)
    \rput*{0}(.7,-.5){ \begin{scshape} Double Cylinder \end{scshape}}
        \endpspicture

         \\

        \multicolumn{2}{c}{\ }\\
        \multicolumn{2}{c}{\ }\\
        \multicolumn{2}{c}{\ }\\
        \multicolumn{2}{c}{\ }\\
        \multicolumn{2}{c}{\ }\\

\end{tabular}

\begin{tabular}{@{\extracolsep{1 cm}} ccc}

    %Slab Lens
    \rput{0}(0.,1.9){
    \epsfysize=.80in
        \epsffile{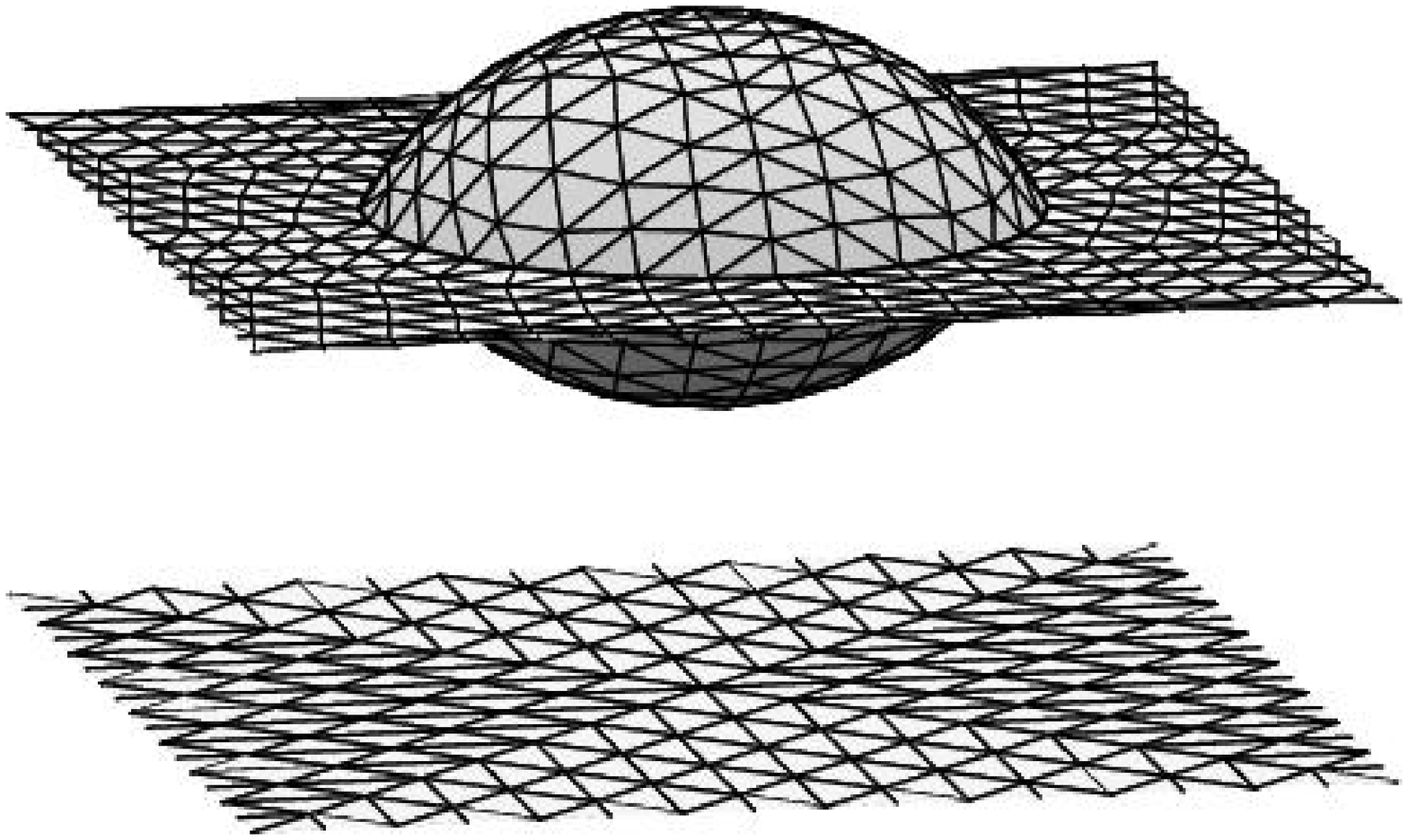}
        }
    \pspicture(-1,-.75)(1,.75)
    \rput*{0}(-1.,-.5){ \begin{scshape} Slab Lens \end{scshape}}
        \endpspicture

        &

        %Center Bubble
    \rput{0}(1.2,1.9){
    \epsfysize=.80in
        \epsffile{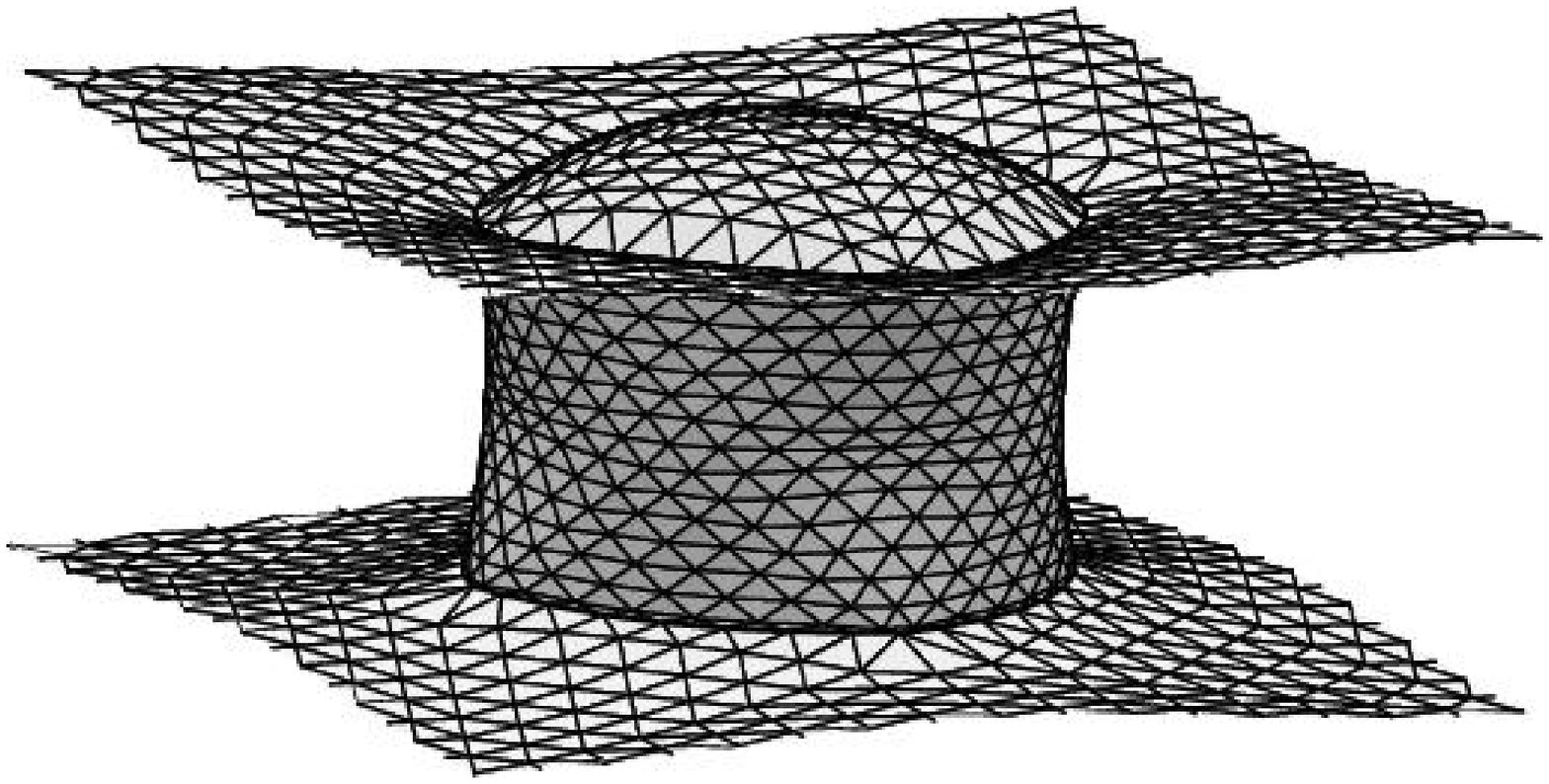}
        }
        \pspicture(-1,-.75)(1,.75)
    \rput*{0}(0.,-.5){ \begin{scshape} Center Bubble \end{scshape}}
        \endpspicture

        &

        %Cylinder String
    \rput{270}(1.8,1.9){
    \epsfysize=1.2in
        \epsffile{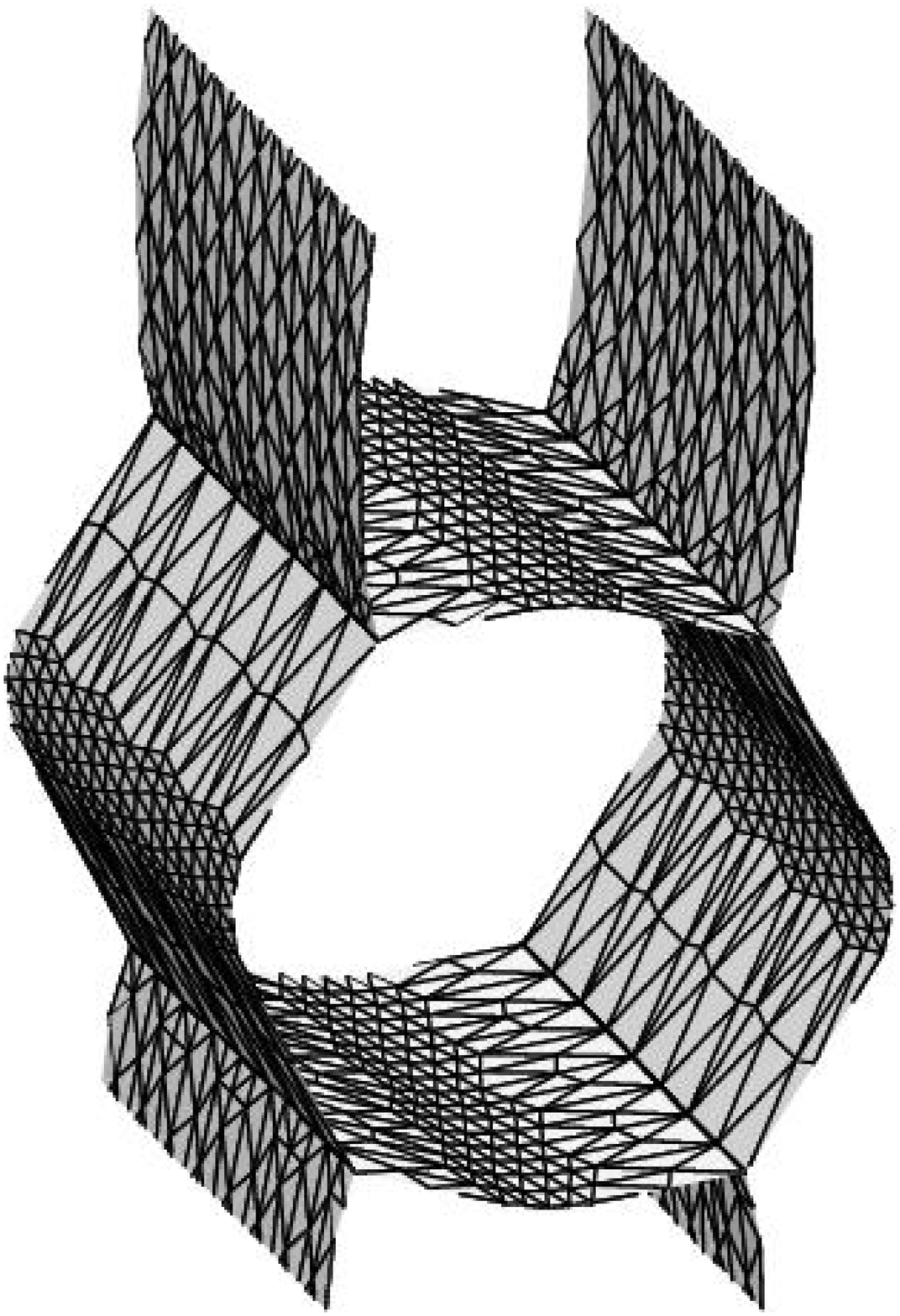}
        }
        \pspicture(-1,-.75)(1,.75)
    \rput*{0}(.7,-.5){ \begin{scshape} Cylinder String \end{scshape}}
        \endpspicture
    \\
        \multicolumn{2}{c}{\ }\\
        \multicolumn{2}{c}{\ }\\
        \multicolumn{2}{c}{\ }\\
        \multicolumn{2}{c}{\ }\\
        \multicolumn{2}{c}{\ }\\

\end{tabular}

\begin{tabular}{@{\extracolsep{.5 cm}} cc}

    %Slab Cylinder
    \rput{-20}(0.2,1.9){
    \epsfysize=1.0in
        \epsffile{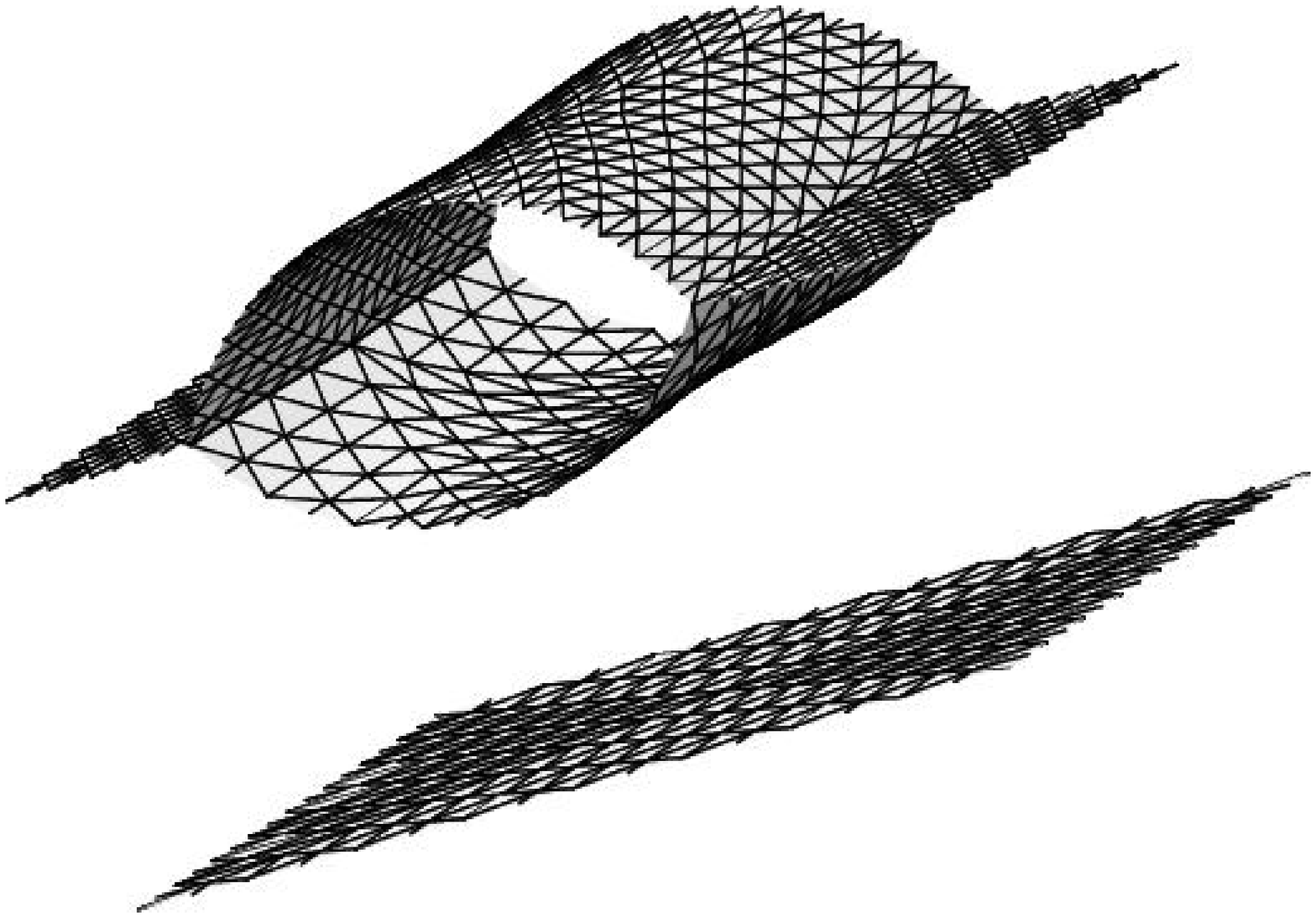}
        }
    \pspicture(-1,-.75)(1,.75)
        \rput*{0}(-.8,-.5){ \begin{scshape} Slab Cylinder \end{scshape}}
        \endpspicture

        &

        %Double Slab
    \rput{-10}(2.,1.9){
    \epsfysize=1.4in
        \epsffile{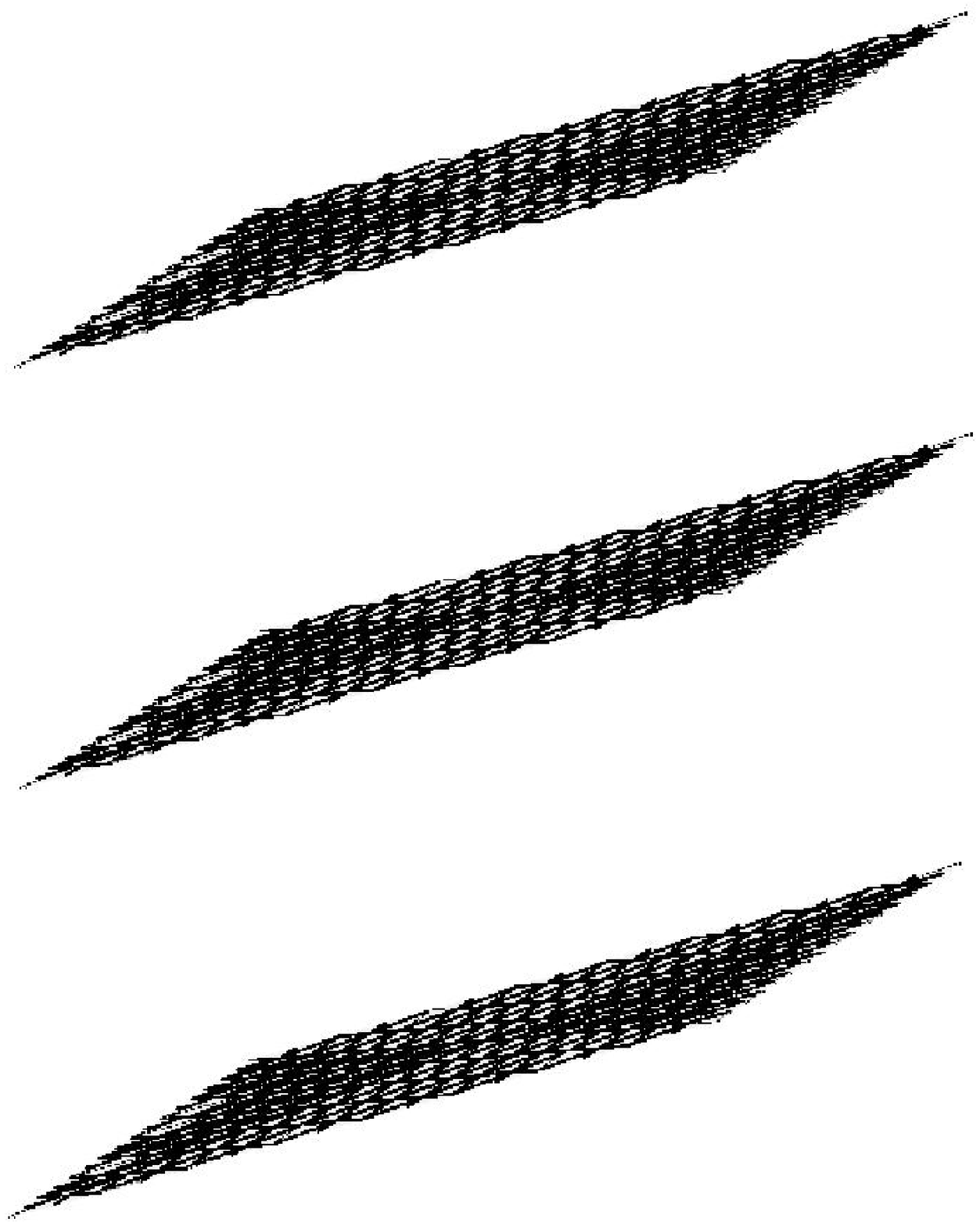}
        }
        \pspicture(-1,-.75)(1,.75)
        \rput*{0}(.7,-.5){ \begin{scshape} Double Slab \end{scshape}}
        \endpspicture

\end{tabular}
\end{center}
\caption{\label{catalog}Catalog of Conjectured Minimizers.}
\end{figure}
%%%%%%%%%%%%%%%%%%%%%%%%%%%%%%%%%%%%%%%%%%%%%%%%%%%%%%%%

The \emph{double bubble problem} is a two-volume generalization of the
famous isoperimetric problem. The isoperimetric problem seeks the
least-area way to enclose a single region of prescribed volume. About
200 BC, Zenodorus argued that a circle is the least-perimeter
enclosure of prescribed area in the plane (see~\cite{Heath}). In 1884,
Schwarz ~\cite{Schwarz} proved by symmetrization that a round sphere
minimizes perimeter for a given volume in \Rth. Isoperimetric problems
arise naturally in many areas of modern mathematics. Ros
~\cite{RosMSRI} provides a beautiful survey.

Soap bubble clusters seek the least-area (least-energy) way to enclose
and separate several given volumes. Bubble clusters have served as
models for engineers, architects, and material scientists (the chapter
by Emmer in \cite{arch_and_math} is a nice survey of architectual
applications, and the text by Weaire and Hutzler \cite{foams}, an
introduction to the physics of space-filling bubble clusters or
\emph{foams}, discusses numerous other applications).  \\

\paragraph{{\bf Existence and regularity }}
The surface area minimizing property of bubble clusters can be
codified mathematically in various useful ways, using the rectifiable
currents, varifolds, or $(\mathbf{M}, \epsilon, \delta)$-minimal sets
of geometric measure theory (see ~\cite{MorganF}). In three
dimensions, mathematically idealized bubble clusters consist of
constant-mean-curvature surfaces meeting smoothly in threes at
$\onetwenty$ along smooth curves, which meet in fours at a fixed angle
of approximately $109^{\circ}$ (\cite[Theorems II.4, IV.5,
IV.8]{Taylor}, or see ~\cite[Section 13.9]{MorganF}).  \\

\paragraph{{\bf Recent results}}
The existence and regularity of solutions to the double bubble problem
played a key role in the proof by Hutchings, Morgan, Ritor\'e, and Ros
(\cite{HMRR}, \cite{HMRRelec}; see \cite[Chapter 14]{MorganF}) that
the standard double bubble, the familiar shape consisting of three
spherical caps meeting one another at 120 degree angles, provides the
least-area way to separate two volumes in \Rth. The proof relied on a
component bound that had been developed by Hutchings \cite[Theorem
4.2]{H}, and a detailed stability argument to rule out the possible
remaining candidates. Reichardt \emph{et al.} ~\cite{rei} extended these
results to \Rf.

Contemporaneously with the development of results presented here,
Corneli, Holt, Leger, and Schoenfeld \cite{CHLS} produced a more or
less complete solution to the \Ttw\ version of the problem. Their
proofs also rely on regularity theory, which in two dimensions implies
that bubbles are bounded by circular arcs. A variational bound on the
number of components of bubble clusters in surfaces due to Wichiramala
\cite{MW} provided considerable added simplifications. A sequence of
proofs using the techniques of plane geometry then eliminated all but
five candidates (of which only four are expected to appear). Examining
the \Ttw\ candidates was helpful to us in our work on the \Tth\
project. \\

\paragraph{{\bf Bubbles in the three-torus }}
In comparison with \Rth\ and \Ttw, the double bubble problem in \Tth\
appears to be more difficult. In the torus it is not possible to push
through a component bound like Hutchings', since a key step in his
proof is to show that the double bubble has an axial symmetry. Nor is
a variational bound after Wichiramala forthcoming, due to the
additional topological complications in three dimensions. Until some
new approach provides a component bound, there probably will be no
definitive results. Indeed, the single bubble for the three-torus is
not yet completely understood, although there are partial results. The
smallest enclosure of half of the volume of the torus was shown by
Barthe and Maurey \cite[Section 3]{Barthe} to be given by two parallel
two-tori. Morgan and Johnson \cite[Theorem 4.4]{MJ} show that the
least-area enclosure of a small volume is a sphere. Spheres, tubes
around geodesics, and pairs of parallel two-tori are shown to be the
only types of area minimizing enclosures for most tori by the work of
Ritor\'e and Ros (\cite[Theorem 4.2]{Rit}, \cite{RitRos}). \\

\paragraph{{\bf Theorem about small volumes }}
Theorem 4.1 states that any sequence of area minimizing double bubbles
of decreasing volume and fixed volume ratio has a tail consisting of
standard double bubbles. The central difficulty is to bound the
curvature. This accomplished, we show that the bubble lies inside some
small ball that lifts to \Rth, where a minimizer is known to be
standard \cite{HMRR}. The result extends to any flat 3- or 4-manifold
with compact quotient by the isometry group.

Our proof goes roughly as follows. From the original sequence of
double bubbles we generate a new sequence by rescaling the manifold at
each stage so that one of the volumes is always equal to one. We can
then apply compactness arguments and area estimates to the rescaled
sequence to show that certain subsequences of translates have
non-trivial limits. These limits are used to obtain a curvature bound
on the original sequence. With such a bound, we can apply monotonicity
to conclude that if the volumes are small, the double bubble is
contained in a small ball. We conclude that it must be the same as the
minimizer in \Rth\ or \Rf, \emph{i.e.} that it must be the standard
double bubble by \cite{HMRR} or \cite{rei}). \\

\paragraph{{\bf Plan of the Paper.}}
Section 2 reviews the methods leading to our Central Conjecture 2.1
and to Figures 1 and 2. Section 3 surveys some subconjectures. Section
4 focuses on the proof of Theorem 1 on small volumes. Section 5 shifts
from the cubic to other tori and discusses other conjectures and
candidates, including a ``Hexagonal Honeycomb." \\

\paragraph{{\bf Acknowledgments.}}
This paper has its origins in a problem given by Frank Morgan to the
audience of his two week course on Geometric Measure Theory and the
Proof of the Double Bubble Conjecture, at the Clay Mathematics
Institute (CMI) Summer School on the Global Theory of Minimal
Surfaces, held at the Mathematical Sciences Research Institute (MSRI)
in June and July of 2001.

We owe our deepest gratitude to Frank Morgan for introducing us to
this problem and taking the time to guide us along as we have worked
on it, by many helpful discussions and constructive comments. CMI and
MSRI are to be thanked for making the Summer School possible, and for
providing research and travel support to the authors. We should like
to extend our thanks to Joel Hass, David Hoffman, Arthur Jaffe,
Antonio Ros, Harold Rosenberg, Richard Schoen, and Michael Wolf, who
organized the Summer School. We benefited tremendously from the
outstanding lectures and informal discussions, for which more people
are owed thanks than we can name here. Three faculty participants,
Michael Dorff, Denise Halverson, and Gary Lawlor, of Brigham Young
University, are to be singled out and thanked for helping us get our
computations started. We thank Manuel Ritor\'e and John Sullivan for
helpful conversations. Other participants who made notable
contributions include Baris Coskunuzer, Brian Dean, Dave Futer, Tom
Fleming, Jim Hoffman, Jon Hofmann, Paul Holt, Matt Kudzin, Jesse
Ratzkin, Eric Schoenfeld, and Jean Steiner.

Yair Kagan of New College of Florida created the picture of the
Hexagonal Honeycomb in Figure 4. The authors would like thank him for
this contribution.

Corneli participated in the 2001 ``SMALL" undergraduate summer
research Geometry Group, supported by Williams College and by the
National Science Foundation through an REU grant to Williams and an
individual grant to Morgan.

%*PHASE DIAGRAM
\begin{figure}[h]
\begin{center}
\includegraphics[width=.95\textwidth]{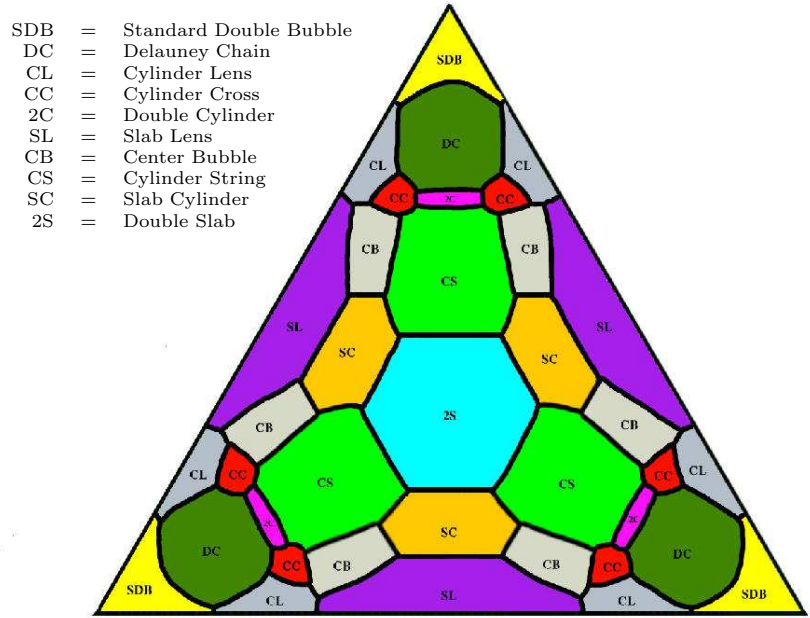}
\rput{0}(-9.5,7)
{
{
\scriptsize
\begin{tabular}{rcl}
SDB &=& Standard Double Bubble\\
DC &=& Delauney Chain\\
CL &=& Cylinder Lens\\
CC &=& Cylinder Cross\\
2C &=& Double Cylinder\\
SL &=& Slab Lens\\
CB &=& Center Bubble\\
CS &=& Cylinder String\\
SC &=& Slab Cylinder\\
2S &=& Double Slab\\
\end{tabular}
}
}
\end{center}
\caption{
\label{phase_portrait}
Phase portrait: volumes and corresponding double bubble. In the
center both regions and the complement have one third of the total
volume; along the edges one volume is small; in the corners two
volumes are small.
}
\end{figure}

\section{The Conjecture}
\paragraph{{\bf Generating Candidates.}}
Many possibilities for double bubbles in the three-torus \Tth\ were
proposed in brainstorming sessions by participants in the Clay/MSRI
Summer School. In order to classify the candidates we used the
following method. Starting with a standard double bubble, we imagined
one of the two volumes growing until the bubble enclosing it wrapped
around the torus and encountered an obstruction. Following the
principals of regularity for bubble clusters, if a bubble collided
with itself we opened the walls up, whereas if two different bubbles
collided we allowed them to stick together. We then repeated this
procedure for these new double bubbles, sometimes changing the
perspective slightly slightly or becoming a bit more fanciful
(e.g. the Center Cylinder of Figure 3 or Gary Lawlor's Fire Hydrant of
Figure 3, also known as ``Scary Gary"). \\

\paragraph{{\bf Producing the Phase Diagram.}}
Brakke's \SE\ ~\cite{SE} was used to closely approximate the minimal
area that a double bubble of each type needs to enclose specified
volumes. Our initial simulations gave us the approximate surface area
for each candidate double bubble, tested on partitions $v_1\colon
v_2\colon v_3$ of a unit volume taken in increments of $0.01$. From
the data obtained in these simulations, we found the least-area
competitor for each volume triple. Figure~1 shows the candidates we
found to be minimizing for some set of volumes. The phase diagram
appearing in Figure \ref{phase_portrait} is the result of refining our
initial computations along the boundaries with a 0.005 increment. \\

\begin{centralconj} \label{central}
The ten double bubbles pictured in Figure 1 represent each type of
surface area minimizing two-volume enclosure in a flat, cubic
three-torus, and these types are minimizing for the volumes
illustrated in Figure 2.
\end{centralconj}

\paragraph{{\bf Comments.}}
One might expect that minimizers would be found among the various
regularity-satisfying conglomerations of topological spheres and
products of spheres and homotopically non-trivial tori, other
possibilities being excessively complex. This was borne out in our
computations. It is interesting to note, however, that not all of the
simple possibilities along appeared as minimizers, for example the
Transverse Cylinders pictured in Figure 3. The various double bubbles
of Figure 3, while stable for a certain range of volumes, are never
area minimizing. A challenging unsolved problem is to find all of the
stable non-minimizing bubbles. Note that a given type from Figure 1
might be stable for a much wider range of volumes than those for which
it actually minimizes surface area.

It is worth observing that all of the conjectured minimizers for the
double bubble problem on \Ttw\ (a Standard Double Bubble, Band with
Lens, Symmetric Chain, and Double Band) are echoed here in at least
two ways. The Double Cylinder, the Slab Cylinder, the Double Slab, and
the Cylinder String are \Ttw\ minimizers $\times$ \Ton. There are
also more direct analogues, as is seen by comparing, for example, the
three- and two-dimensional Standard Double Bubbles, or the D\'elauney
and Symmetric Chains. See Corneli \emph{et al.} \cite{CHLS} for more
on the $T^2$ minimizers.

\begin{figure}[h]
\begin{center}
    \pspicture(-1,-3.9)(1,3.7)
 \rput{0}(1.1,0){
\begin{tabular}{@{\extracolsep{2.3 cm}} cc}

    %Transverse Cylinders
    \rput{0}(0.,0){
    \epsfysize=.80in
        \epsffile{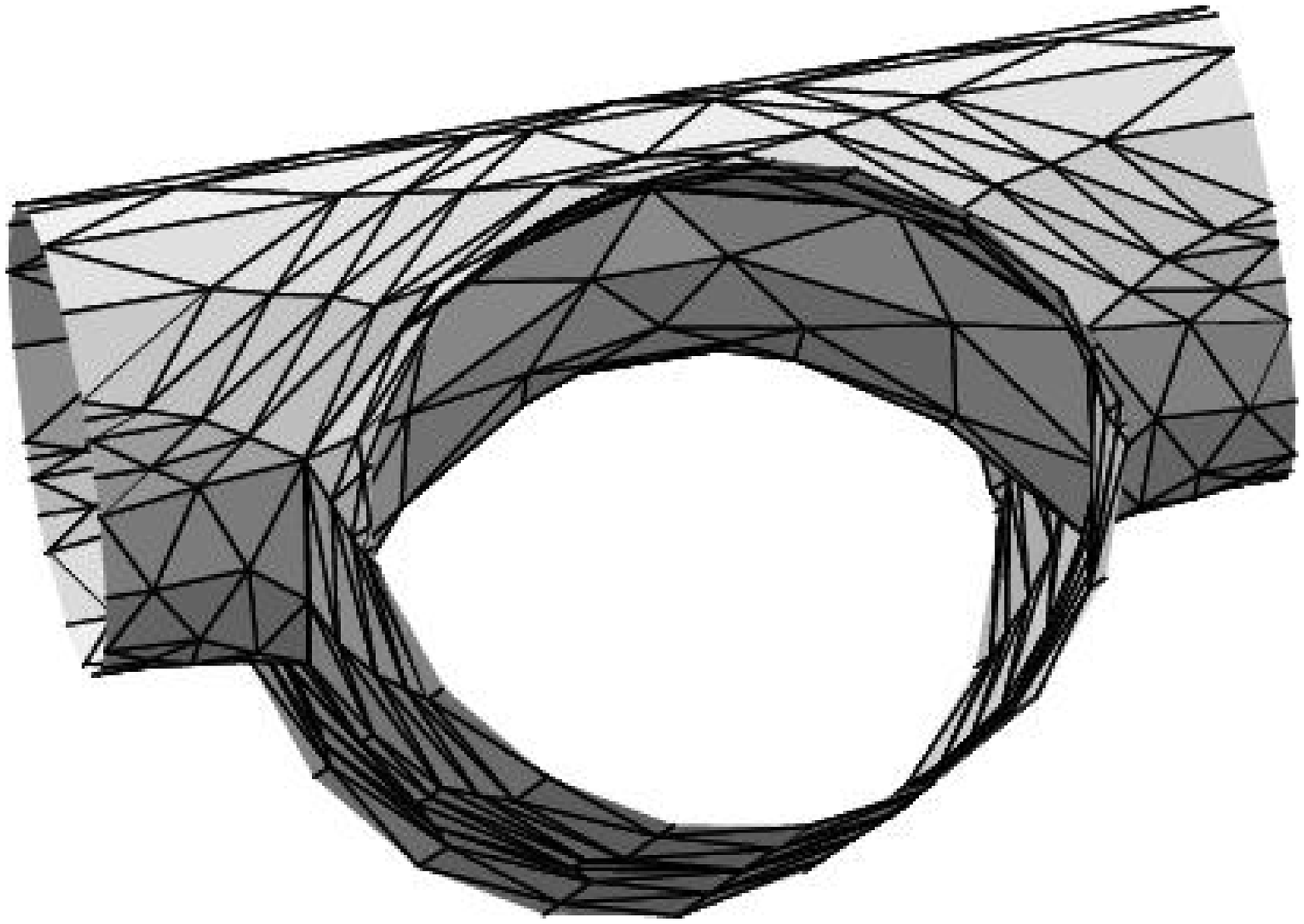}
        }
    \pspicture(-1,-.75)(1,.75)
    \rput*{0}(-1.,-2.7){ \begin{scshape} Transverse Cylinders \end{scshape}}
        \endpspicture

        &

        %Stump
    \rput{0}(0,0){
    \epsfysize=.80in
        \epsffile{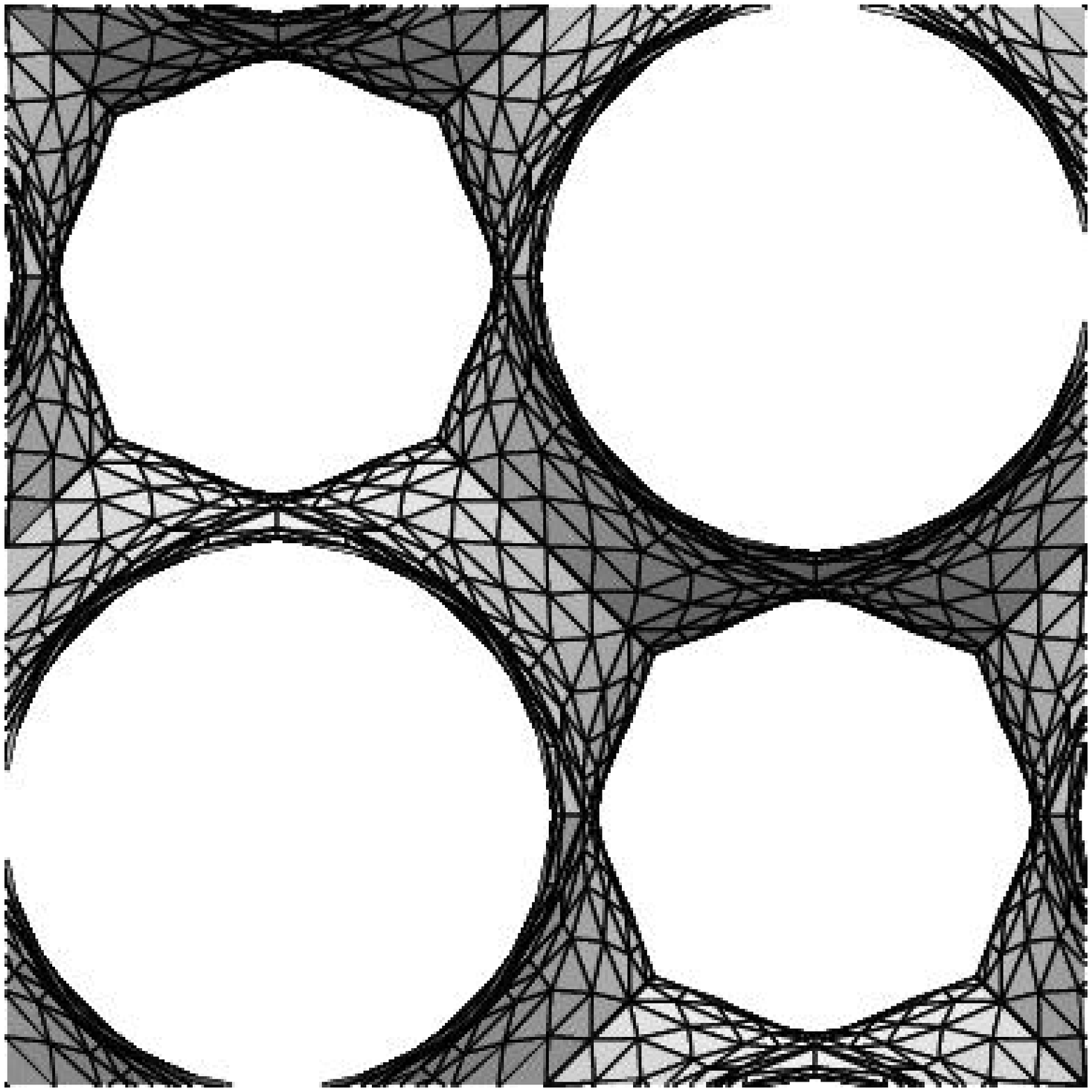}
        }
        \pspicture(-1,-.75)(1,.75)
    \rput*{0}(-1,-2.7){ \begin{scshape} Double Hydrant \end{scshape}}
        \endpspicture   \\

         \multicolumn{2}{c}{\ }\\
        \multicolumn{2}{c}{\ }\\
         \multicolumn{2}{c}{\ }\\
        \multicolumn{2}{c}{\ }\\
        \multicolumn{2}{c}{\ }\\
        \multicolumn{2}{c}{\ }\\

        %Center Cylinder
    \rput{0}(0,0){
    \epsfysize=1.2in
        \epsffile{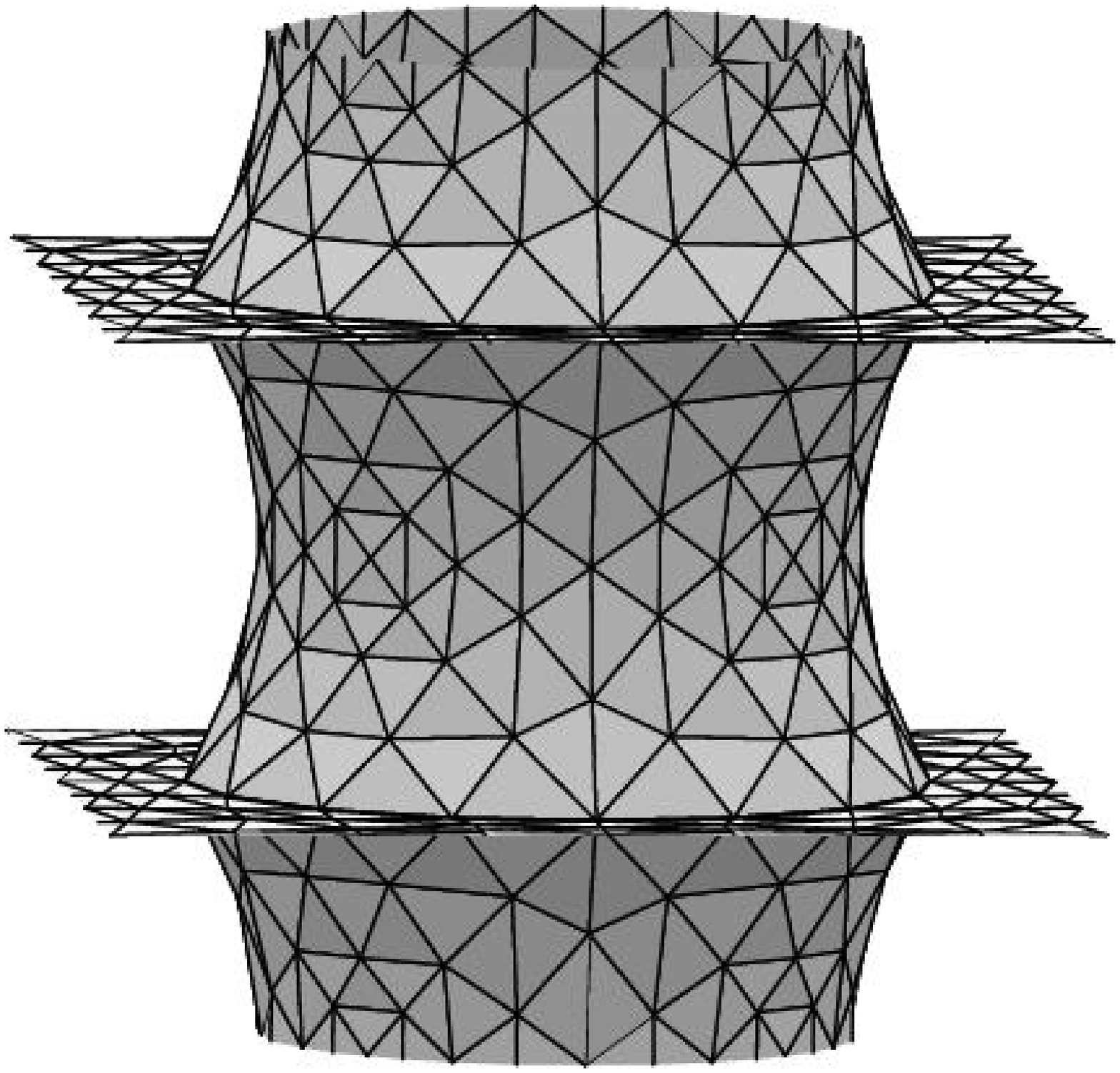}
        }
        \pspicture(-1,-.75)(1,.75)
    \rput*{0}(-1.1,-2.9){ \begin{scshape} Center Cylinder \end{scshape}}
        \endpspicture
&
        %scary gary
    \rput{0}(0,0){
    \epsfysize=1.2in
        \epsffile{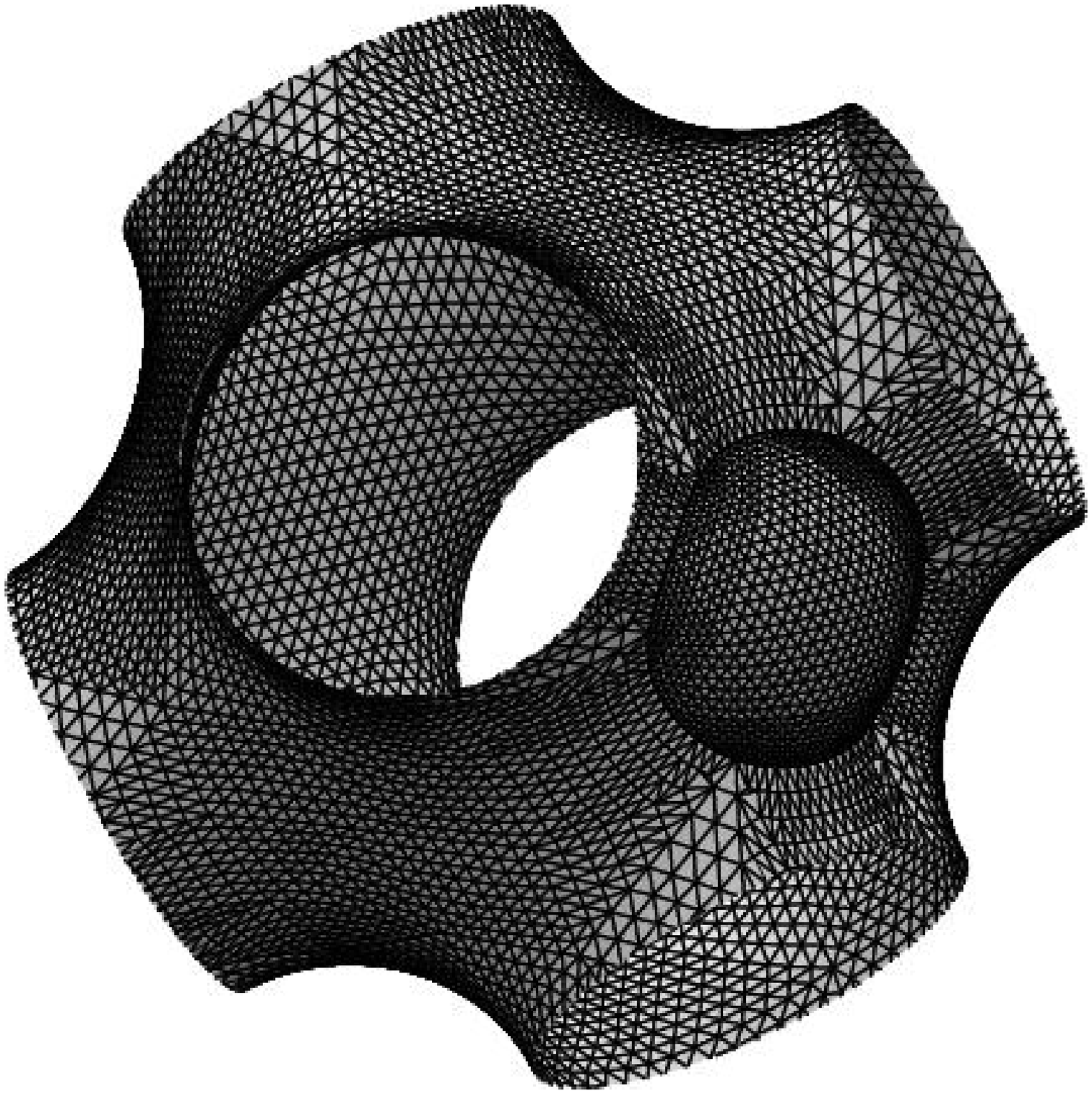}
        }
        \pspicture(-1,-.75)(1,.75)
    \rput*{0}(-1,-2.9){ \begin{scshape} Hydrant Lens \end{scshape}}
        \endpspicture

    \\

        \multicolumn{2}{c}{\ }\\
        \multicolumn{2}{c}{\ }\\
        \multicolumn{2}{c}{\ }\\
        \multicolumn{2}{c}{\ }\\
        \multicolumn{2}{c}{\ }\\

\end{tabular}
}   
\endpspicture

\end{center}
\caption{Inefficient double bubbles}
\end{figure}

\section{Subconjectures}
We now present a list of natural subconjectures, suggested either by
the phase diagram or by examination of the pictures made
with \SE.

One immediate observation is that the edges of our phase diagram
appear to characterize single bubbles in \Tth.

\begin{conj}[Ritor\'e and Ros 
\cite{Rit}, \cite{RitRos}, \cite{RosMSRI}] The optima for the
isoperimetric problem in a cubic \Tth\ are the sphere, cylinder, and
slab.
\end{conj}

We still do not know proofs for the following intuitive conjectures
about the double bubble problem (although Theorem 4.1 gives partial
results on Conjecture 3.3):

\begin{conj} \label{connectedness}
An area-minimizing double bubble in \Tth\ has connected regions and
complement.
\end{conj}

\begin{conj} \label{conj:small_volumes}
For two small volumes the standard double bubble is optimal.
\end{conj}

\begin{conj}
For one very small volume and two moderate volumes, the Slab Lens is
optimal.
\end{conj}

\begin{conj}
The first phase transition as small equal, or close to equal, volumes
grow is from the standard double bubble to a chain of two bubbles
bounded by D\'elaunay surfaces.
\end{conj}

D\'elaunay surfaces are constant-curvature surfaces of rotation, and
as such have full rotational symmetry. From the \SE\ pictures, it
appears that the conjectured surface area minimizers always have the
maximial symmetry, given the constraints , a fact which leads to the
following natural conjecture.

\begin{conj} 
A minimizer is as symmetric as possible, given its topological type.
\end{conj}

We will conclude this section with a proposition which represents the
first step towards establishing a symmetry property. Specifically, we
prove that for any double bubble, there is a pair of parallel planes
(actually two-tori) that cut both regions in half. Hutchings
\cite[Theorem 2.6]{H} used the fact that in \Rn\ a double bubble has
two perpendicular planes that divide both regions in half. This is an
easy step in a difficult proof that the function that gives the
least-area to enclose two given volumes is concave. Before proving a
basic result of a similarly elementary flavor for the torus, we
mention that we conjecture that the much deeper concavity result also
holds:

\begin{conj}
The least area to enclose and separate two given volumes in the
three-torus is a concave function of the volumes.
\end{conj}

If one could prove Conjecture 3.7, one would then be able to apply
other ideas in Hutchings' paper \cite[Section 4]{H} to obtain a
functional bound on the number of components of a minimizing bubble.

\begin{prop} [Deluxe Ham Sandwich Theorem\footnote{Name suggested by Eric Schoenfeld, who gave an
independent proof for $T^2$.}] In a rectangular torus, if a double
bubble lies inside a cylinder $S^1 \times D^2$, then there is a plane
that cuts both volumes in half.
\end{prop}

\begin{proof}
This is a generalization of the standard argument for the ham sandwich
theorem in Euclidean space. Assume that the identification occurs in the
vertical direction. Take a circle in the $xy$ plane indexed by $\alpha \in
[0, \pi]$. Then for each $\alpha$, rotate the surface by $\alpha$ in the
$xy$ plane and consider the family of pairs of flat tori parallel to the
$xz$ plane that are at distance 1/2 from each other. Then there is at
least one such pair of planes that cuts the volume $V_1$ in half, for each
$\alpha$. There may be an interval of such pairs for a given $\alpha$.
However, some one parameter family of these planes may be chosen which 
varies continuously as a function of $\alpha$, and one such plane cuts
$V_2$ in half.
\end{proof}

\section{Small Volumes} \label{S:small_volumes}

Conjecture \ref{conj:small_volumes} stated that small volumes are best
enclosed by a standard double bubble. Theorem 4.1 proves that for any
fixed volume ratio, the standard double becomes optimal when the
volumes are sufficiently small. This result holds for many flat three-
and four-dimensional manifolds (and see Remark \ref{in_general}).

The standard double bubble consists of three spherical caps meeting at
120 degrees. (If the volumes are equal, the middle surface is planar.)
It is known to be minimizing in \Rth\ \cite{HMRR} and \Rf\ \cite{rei}.

\begin{thm} \label{thm:small_volumes}
Let $M$ be a smooth flat Riemannian manifold of dimension three or
four, such that $M$ has compact quotient by the isometry group. Fix
$\lambda \in (0,1]$. Then there is an $\epsilon >0$, such that if $0<
v < \epsilon$, a perimeter-minimizing double bubble of volumes $v,
\lambda v$ is standard.
\end{thm}

\begin{rem}
Solutions to the double bubble problem exist for all volume pairs
in manifolds with compact quotient by their isometry group. The
proof is the same as in Morgan ~\cite[Section 13.7]{MorganF}. 
\end{rem}

For the proof we will regard a double bubble as a pair of $3$- or
$4-$dimensional rectifiable currents, $R_1$ and $R_2$, each of
multiplicity one, of volumes $V_1 = \mathbf{M}(R_1)$ and $V_2 =
\mathbf{M}(R_2)$. The total perimeter of such a double bubble is
$\frac{1}{2}(\mathbf{M}(\partial R_1) +\mathbf{M}(\partial R_2)
+\mathbf{M}(\partial(R_1 +R_2)))$. Here $\mathbf{M}$ denotes the mass
of the current, which can be thought of as the Hausdorff measure of
the associated rectifiable set (counting multiplicities). For a review
of the pertinent definitions, see notes from Morgan's course at the
Summer School \cite{MorganCourse} or the texts by Morgan
\cite{MorganF} or Federer \cite{Federer}.

Before proceeding it is helpful to fix some notation. By the Nash and
subsequent embedding theorems, we may assume $M$ is a submanifold of
some fixed \RN. We will consider a sequence of perimeter-minimizing
double bubbles in $M$, containing the volumes $v$ and $\lambda v$, as
$v \rightarrow 0$. For each $v$, $M_v$ will denote $s_v(M)$ in
$s_v(\mRN)$, where $s_v$ is the scaling map that takes regions with
volume $v$ to similar regions with volume $1$. In particular, $s_v$
maps our perimeter-minimizing double bubble containing volumes $v$ and
$\lambda v$ to a double bubble which we call $S_v$ which contains
volumes $1$ and $\lambda$, and which is of course perimeter-minimizing
for these volumes.

We focus on the case of dimension three; the proof for dimension four
is essentially identical.

\begin{lem} 
There is a $\gamma >0 $ such that if $R$ is a region in an open
Euclidean 3-cube $K$ and $\vol(R) \leq \vol(K)/2$, then
%%%%%%%%%%%%%%%%%%%%%%%%%%%%%%%%%%%%%%%%%%%%%%%%%%%%%%%%%%%%
$$\area(\partial R) \geq \gamma  (\vol(R))^{2/3}.$$
%%%%%%%%%%%%%%%%%%%%%%%%%%%%%%%%%%%%%%%%%%%%%%%%%%%%%%%%%%%%
\end{lem}

\begin{proof}
Let $\gamma_0$ be such an isoperimetric constant for a cubical
3-torus, so that $\area(\partial P) \geq \gamma_0 (\vol(P))^{2/3}$ for
all regions $P \subseteq \mTth$. Such a $\gamma_0$ exists by the
isoperimetric inequality for compact manifolds ~\cite[Section
12.3]{MorganF}. Make the necessary reflections and identifications of
the cube $K$ to obtain a torus containing a region $R'$ with eight
times the volume and eight times the boundary area of $R$. The claim
follows, with $\gamma = \gamma_0/2$.
\end{proof}

\begin{proof}[{\bf Proof of Theorem 4.1}]
The first step is to show that the sequence $S_v$, suitably translated
and rotated, has a subsequence that converges as $v \rightarrow 0$ and
has $V_1 \neq 0$ in the limit. Our argument also shows that there is a
subsequence that converges to a limit with $V_2 \neq 0$, but does not
show that there is a subsequence where both volumes are non-zero in
the limit. This is because while we are exerting ourselves trapping
the first volume in a ball, the second one may wander off to infinity.

We first show the existence of a covering $\mathscr{K}_v$ of $ M_v$
with bounded multiplicity, consisting of cubes contained in $M_v$,
each of side-length $L$. Lemma 1 will give us a positive lower bound
on the volume of the part of $R_1$ that is inside one of these cubes
for each $S_v$. We will then apply a standard compactness theorem to
show that a subsequence of the $S_v$'s converges.

Take a maximal packing by balls of radius $\frac{1}{4}
L$. Enlargements of radius $\frac{1}{2} L$ cover $M_v$. Circumscribed
cubes of edge-length $L$ provide the desired covering $\mathscr{K}_v$.
To see that the multiplicity of this covering is bounded, consider a
point $p \in M_v$. The ball centered at $p$ with radius 2$L$ contains
all the cubes that might cover it, and the number of balls of radius
$\frac{1}{4} L$ that can pack into this ball is bounded, implying that
the multiplicity of $\mathscr{K}_v$ is also bounded by some $m>0$.

Now consider some such covering, with $L=2$. By Lemma 1, there is an
isoperimetric constant $\gamma$ such that
%%%%%%%%%%%%%%%%%%%%%%%%%%%%%%%%%%%%%%%%%%%%%%%%%%%%%%%%%%%
$$ \area(\partial (R_1 \cap K_i)) \geq \gamma (\vol(R_1 \cap
K_i))^{2/3}, $$
%%%%%%%%%%%%%%%%%%%%%%%%%%%%%%%%%%%%%%%%%%%%%%%%%%%%%%%%%%%
and therefore, since $\max_k \vol (R_1 \cap K_k) \geq \vol(R_1 \cap
K_i)$ for any $i$,
%%%%%%%%%%%%%%%%%%%%%%%%%%%%%%%%%%%%%%%%%%%%%%%%%%%%%%%%%%%
\begin{equation}
\area (\partial (R_1 \cap K_i) ) \geq \gamma \frac{\vol(R_1 \cap
K_i)}{(\max_k \vol (R_1 \cap K_k) )^{1/3}}. \label{sum_me}
\end{equation}
%%%%%%%%%%%%%%%%%%%%%%%%%%%%%%%%%%%%%%%%%%%%%%%%%%%%%%%%%%%

Note that the total area of the surface is greater than $1/m$ times
the sum of the areas in each cube, and the total volume enclosed is
less than the sum of the volumes, so summing Equation \ref{sum_me}
over all the cubes $K_i$ in the covering $\mathscr{K}_v$ yields
%%%%%%%%%%%%%%%%%%%%%%%%%%%%%%%%%%%%%%%%%%%%%%%%%%%%%%%%%%%
$$\area(S_v) \geq \area(\partial R_1) \geq m \gamma \frac{V_1}{(\max_k
\vol (R_1 \cap K_k))^{1/3}} $$
%%%%%%%%%%%%%%%%%%%%%%%%%%%%%%%%%%%%%%%%%%%%%%%%%%%%%%%%%%%
and
%%%%%%%%%%%%%%%%%%%%%%%%%%%%%%%%%%%%%%%%%%%%%%%%%%%%%%%%%%%
$$(\max_k(\vol(R_1 \cap K_k)))^{1/3} \geq m \gamma
\frac{V_1}{\area(S_v)} \geq \delta >0,$$
%%%%%%%%%%%%%%%%%%%%%%%%%%%%%%%%%%%%%%%%%%%%%%%%%%%%%%%%%%%
because $V_1 = 1$ and it is easy to show that $\area(S_v)$ is bounded
(since there is a bounded way of enclosing the volumes). Translate
each $M_v$ so that the cube where the maximum occurs is centered at
the origin of \RN, and rotate so that the tangent space of each $M_v$
at the origin is equal to a fixed \Rth\ in \RN. The limit of the $M_v$
will be equal to this \Rth. Since a cube with edge-length $L$ centered
at the origin fits inside a ball of radius $2L$ centered at the
origin, we have
%%%%%%%%%%%%%%%%%%%%%%%%%%%%%%%%%%%%%%%%%%%%%%%%%%%%%%%%%%%
$$ \vol(R_1 \cap B(0,2L)) \geq \delta^3 $$
%%%%%%%%%%%%%%%%%%%%%%%%%%%%%%%%%%%%%%%%%%%%%%%%%%%%%%%%%%%
for every $S_v$. By the compactness theorem for locally integral
currents (\cite[pp. 64,88]{MorganF},~\cite[Section 27.3, 31.2,
31.3]{Simon}) we know that a subsequence of the $S_v$ has a limit,
which we will call $D$, with the property that $\vol(R_1) \geq
\delta^3$. This completes the first step.

Since $D$ is contained in the limit of the $M_v$, namely, the copy of
\Rth\ chosen above, and each $S_v$ is minimizing for its volumes, a
standard argument shows that the limit $D$ is the perimeter-minimizing
way to enclose and separate the given volumes in \Rth\
(cf. \cite[13.7]{MorganF}). In the limit, $V_2$ could be zero, in
which case D is a round sphere. If both volumes are non-zero, D is the
standard double bubble (\cite{HMRR} and \cite{rei}, or see
\cite[Chapter 14]{MorganF}).

Our goal is now to prove that the double bubble of volumes $v$,
$\lambda v$ lies inside a trivial ball, when $v$ is small enough.
This can be accomplished using the monotonicity theorem for mass ratio
~\cite[Section 5.1(1)]{Allard}, which implies that for a
perimeter-minimizing bubble cluster, a small ball around any point on
the surface contains some substantial amount of area. This will limit
the number of disjoint balls we can place on the surface. The
monotonicity theorem applies only to surfaces for which the mean
curvature is bounded, \emph{i.e.}, for which there is a $C$ such that
for smooth variations
%%%%%%%%%%%%%%%%%%%%%%%%%%%%%%%%%%%%%%%%%%%%%%%%%%%%%%%%%%%%
$$ \frac{dA}{dV} \leq C.$$
%%%%%%%%%%%%%%%%%%%%%%%%%%%%%%%%%%%%%%%%%%%%%%%%%%%%%%%%%%%%

Accordingly, the second major step in the proof is to obtain such
bound on the curvature, for $v$ small. It suffices to show that all
smooth variation vector fields have the property that changes in the
volume of $S_v$ and in the area of $\partial S_v$ are controlled. Take
a smooth variation vector field $F$ in \Rn\ such that for $D$,
%%%%%%%%%%%%%%%%%%%%%%%%%%%%%%%%%%%%%%%%%%%%%%%%%%%%%%%%%%%%
$$ dV_1/dt = \int_{\partial R1} (F \cdot n)\ dA= c > 0 $$
%%%%%%%%%%%%%%%%%%%%%%%%%%%%%%%%%%%%%%%%%%%%%%%%%%%%%%%%%%%%
and
%%%%%%%%%%%%%%%%%%%%%%%%%%%%%%%%%%%%%%%%%%%%%%%%%%%%%%%%%%%%
$$ dV_2/dt = 0.$$
%%%%%%%%%%%%%%%%%%%%%%%%%%%%%%%%%%%%%%%%%%%%%%%%%%%%%%%%%%%%
(Note that we need the first volume to be non-zero, or its variation
could be zero.) For $v$ small enough the subsequence of $S_v$ headed
towards the limit $D$ has the property that $dV_1/dt$ is approximately
$c$ and $dV_2/dt$ is approximately 0.

By the argument in the first step, we can translate each $S_v$
similarly, so that a subsequence of the subsequence above converges to
a minimizer $D'$ in \Rth\ where the second volume is non-trivial. This
time we take a smooth variation vector field $F'$, such that for $v$
small enough the subsequence of $S_v$ headed to $D'$ has the property
that $dV_1/dt$ is approximately 0 and $dV_2/dt$ is approximately $c' >
0$. This proves that for this subsequence the change in volume is
bounded below.

Now we need to show that the change in area is bounded above. This
follows from the fact that every rectifiable set can be thought of as
a varifold ~\cite[Section 11.2]{MorganF}. By compactness for varifolds
~\cite[Section 6]{Allard}, the $ S_v$, situated so that the first
volume does not disappear, converge as varifolds to some varifold,
$J$. The first variation of the varifolds also converge, \emph{i.e.},
$\delta S_v \rightarrow \delta J$, see ~\cite{Allard}. The first
variation of a varifold is a function representing the change in
area. Therefore, far enough out in the sequence the change in area of
the $S_v$ under $F$ is bounded close to the change in area of $J$
under $F$, which is finite. Similarly, the change in area of the
$S_v$ under $F'$ is bounded.

We conclude that the mean curvature of $S_v$ is bounded for two
independent directions in the two-dimensional space of volume
variations, and hence for all variations. This completes the second
step.

The third and final step is to show that all of the surface area is
contained in some ball in $M_v$, of fixed radius for all
$v$. Eventually, as $v$ shrinks and $M_v$ grows, this ball will have
to be trivial in $M_v$. We will then use the result that the optimal
double bubble in \Rth\ is standard to show that our double bubble is
standard as well.

As the $S_v$ are approaching the minimizer in \Rth, there must be a
bound $A$ on the perimeters of the $S_v$. By monotonicity of mass
ratio ~\cite[Section 5.1(1)]{Allard}, every unit ball centered at a
point of $S_v$ contains perimeter $\delta >0$. Therefore there are at
most $A/\delta$ such disjoint balls.
 
We claim that the $S_v$ are eventually connected. There is an upper
bound on the diameter of any component, $2A/\delta$. Since we are
controlling curvature, our components cannot become too small. Since
every unit ball contains at least $\delta$ area, we also have a lower
bound on the area of each component, when unit balls are
trivial. Unless eventually the $S_v$ are connected, you can arrange to
get in the limit a disconnected minimizer in \Rth, a contradiction.

Hence $S_v$ is contained in a ball of radius $2 A/ \delta$ for all
$v$. Since our original manifold has compact quotient by its isometry
group, there is a radius such that balls in the original manifold of
that radius or smaller are topologically trivial. Hence, as we expand
the manifold, eventually balls of radius $2 A/ \delta$ can be lifted
to \Rth, which means that they are Euclidean. Hence, $S_v$ is
eventually contained in a Euclidean ball, and is therefore the
standard double bubble (\cite{rei}, ~\cite{HMRR}, ~\cite[Chapter
14]{MorganF}).

Finally, since $S_v \subset M_v$ is simply a scaled version of the
original double bubble in $M$, we conclude that the original double
bubble is standard as desired.
\end{proof}

\begin{rem} \label{in_general}
Given $n$, $m$, similar arguments show that for any smooth
$n$-dimensional Riemannian manifold with compact quotient by the
isometry group, given $0 \leq \lambda \leq 1$, there are $C$,
$\epsilon >0$, such that for any $0 < v <\epsilon$, a minimizing
cluster with $m$ prescribed volumes between $\lambda v $ and $v$ lies
inside a ball of diameter at most $C v^{1/n}$.

To further deduce that the cluster smoothly approximates a Euclidean
minimizer would require knowing that convergence weakly and in
measure, under bounded mean curvature, implies $C^1$ convergence, as
is known for hypersurfaces without singularities ([1, Section 8], see
\cite[Section 1.2]{Smallandcompact})
\end{rem}

\section{Special Tori}

Changing the shape of the torus, by stretching it or by skewing some
or all of its angles, would certainly change the phase diagram of
Figure 2:

\begin{conj}
In the special case of a very long $T^3$ the Double Slab is optimal
for most volumes.
\end{conj}

Special tori may have special minimizers:

\begin{conj} 
For the special case of a torus based on a relatively short
$\sixty$-rhombic right prism, the Hexagonal Honeycomb prism of Figure
4 is a perimeter-minimizing double bubble for which both regions and
the exterior each have one third the volume, or when two volumes are
equal and the third is close.
\end{conj}

Indeed, for such volumes the Hexagonal Honeycomb ties the Double Slab,
just as in the $\sixty$-rhombic two-torus a Hexagonal Tiling ties the
Double Band \cite{CHLS}.

\begin{figure}[hbtp]
\centering
\rotatebox{270}{
\includegraphics[bb= 53 58 437 599, width=.5\textwidth, clip]{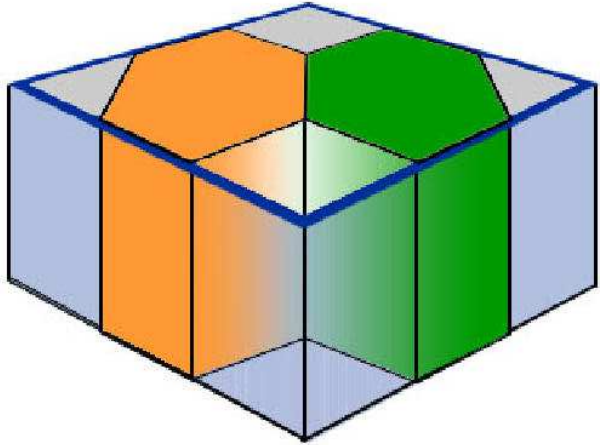}
}
\caption{
\label{fig:HexHoney} Hexagonal Honeycomb}
\end{figure}

In the \emph{triple bubble problem} for the Face Centered Cubic (FCC)
and Body Centered Cubic (BCC) tori we would expect to find minimizers
that lift to \Rth\ as periodic foams with cells of finite volume. The
Weaire-Phelan foam, a counterexample to Kelvin's conjectured best way
to divide space into unit volumes (\cite{WP}, \cite{KS}; see
\cite[Chapter 15]{MorganF}, \cite{BrakkeWP}), would be expected to
appear as lifts of some solutions to the triple bubble problem in the
BCC torus. Kelvin's foam might appear as a solution to the triple
bubble problem in the FCC torus.  One might also see \emph{eight}
Kelvin cells in the BCC torus, or \emph{sixteen} Kelvin cells in the
standard cubic torus.  Eight Weaire-Phelan cells also fit in the
standard cubic torus. Furthermore, by scaling, higher numbers of cells
fit into these tori.  

In contrast with the triple bubble problem and problems with more
volume constraints, it is extremely unlikely that a solution to the
double bubble problem in any torus would lift to a division of \Rth\
into finite volumes, since by regularity, singular curves meet in
fours. This means that singular points look locally like the cone over
a tetrahedral frame, and hence have four volumes coming
together. Since a region will never be adjacent to another component
of the same region (because the dividing wall could be removed to
decrease area and maintain volumes), a foam generated using a
fundamental domain coming from a double double in the three-torus
would have to exhibit the strange property that the singular curves
never meet.

\begin{conj}
There are no least-area divisions of \Tth\ into three volumes that
lift to a foam in \Rth.
\end{conj}

This conjecture suggests that it is not likely that there 
will be are any other special minimizers for the double bubble problem.

\begin{conj}
The double bubbles of Figure \ref{catalog} together with the Hexagonal
Honeycomb of Figure ~\ref{fig:HexHoney} comprise the complete set of area
minimizing double bubbles for all three tori.
\end{conj}

As a final remark, in light of the fact that the triple bubble problem
in the torus seems likely to produce so many interesting candidates,
we would like to mention one final conjecture.

\begin{conj}
For the triple bubble problem in a cubic \Tth\ in the case where one
of the volumes is small, the minimizers will look like the double
bubbles of Figure 2 with a small ball attached. The phase diagram will
look just like our Figure 3.
\end{conj}

\bibliography{mybib.bib}

\begin{thebibliography}{10}

\bibitem{Allard}
William~K. Allard.
\newblock On the first variation of a varifold.
\newblock {\em Ann. of Math.(2)}, 95:418--446, 1972.

\bibitem{Barthe}
F.~Barthe and B.~Maurey.
\newblock Some remarks on isoperimetry of {G}aussian type.
\newblock {\em Ann. Inst. H. Poincar\'e Probab. Statist.}, 36:419--434, 2000.

\bibitem{SE}
K.~Brakke.
\newblock Surface evolver.
\newblock http://www.susqu.edu/facstaff/b/brakke/evolver/evolver.html.

\bibitem{BrakkeWP}
K.~Brakke.
\newblock Century-old soap bubble problem solved!
\newblock {\em Imagine That!}, 3:1--3, Fall, 1993.
\newblock Publication of The Geometry Center.

\bibitem{CHLS}
Joseph Corneli, Paul Holt, Nicholas Leger, and Eric Schoenfeld.
\newblock The double bubble conjecture on the flat 2-torus, 2001.
\newblock Williams College NSF "SMALL" Geometry Group undergraduate research
  report, {\tt http://lanfiles.williams.edu/\~{}fmorgan/geom.group/torus.ps}.

\bibitem{arch_and_math}
Michele Emmer.
\newblock Architecture and mathematics: Soap bubbles and soap films.
\newblock In Kim Williams, editor, {\em Nexus: Architecture and Mathematics},
  pages 53--65. Edizioni dell'Erba, Fucecchio, 1996.

\bibitem{Federer}
Herbert Federer.
\newblock {\em Geometric Measure Theory}.
\newblock Springer-Verlag, Berlin, 1969.

\bibitem{Heath}
Thomas Heath.
\newblock {\em A History of Greek Mathematics}.
\newblock Oxford University Press, Oxford, 1960.
\newblock (Vol. II).

\bibitem{H}
Michael Hutchings.
\newblock The structure of area-minimizing double bubbles.
\newblock {\em J. Geom. Anal.}, 7:285--304, 1997.

\bibitem{HMRRelec}
Michael Hutchings, Frank Morgan, Manuel Ritore, and Antonio Ros.
\newblock Proof of the double bubble conjecture.
\newblock {\em Electron. Res. Announc. Amer. Math. Soc.}, 6:45--49 (elecronic),
  2000.

\bibitem{HMRR}
Michael Hutchings, Frank Morgan, Manuel Ritor{\'e}, and Antonio Ros.
\newblock Proof of the double bubble conjecture.
\newblock {\em Ann. Math. 155}, pages 459--489, March 2002.

\bibitem{KS}
R.~B. Kusner and J.~M. Sullivan.
\newblock Comparing the {W}eaire-{P}helan equal-volume foam to {K}elvin's foam.
\newblock {\em Forma}, 11:233--242, 1996.

\bibitem{MorganF}
Frank Morgan.
\newblock {\em Geometric Measure Theory: a Beginner's Guide}.
\newblock Academic Press Inc., San Diego, CA, third edition, 2000.

\bibitem{Smallandcompact}
Frank Morgan.
\newblock Small perimeter-minimizing double bubbles in compact surfaces are
  standard.
\newblock {\em Electronic Proceedings of the 78th annual meeting of the
  Lousiana/Mississippi Section of the MAA, Univ. of Miss., March 23-24}, 2001.
\newblock to appear.

\bibitem{MJ}
Frank Morgan and David~L. Johnson.
\newblock Some sharp isoperimetric theorems for {R}iemannian manifolds.
\newblock {\em Indiana U. Math J.}, 49:1017--1041, 2000.

\bibitem{MorganCourse}
Frank Morgan and Manuel Ritor\'e.
\newblock Geometric measure theory and the proof of the double bubble
  conjecture.
\newblock Lecture notes by Ritore on Morgan's lecture series at the Clay
  Mathematics Institute Summer School on the Global Theory of Minimal Surfaces,
  at the Mathematical Sciences Research Institute, Berkeley, California 2001
  http://www.ugr.es/~ritore/preprints/course.pdf.

\bibitem{MW}
Frank Morgan and Wacharin Wichiramala.
\newblock The standard double bubble is the unique stable double bubble in
  $\mathbf{R}^2$.
\newblock {\em Proc. AMS, to appear}.

\bibitem{rei}
Ben~W. Reichardt, Cory Heilmann, Yuan~Y. Lai, and Anita Spielmann.
\newblock Proof of the double bubble conjecture in ${R^4}$ and certain higher
  dimensional cases.
\newblock {\em Pacific J. Math., {\textrm to appear}}.

\bibitem{Rit}
Manuel Ritor{\'e}.
\newblock Applications of compactness results for harmonic maps to stable
  constant mean curvature surfaces.
\newblock {\em Math. Z.}, 226:465--481, 1997.

\bibitem{RitRos}
Manuel Ritor\'e and Antonio Ros.
\newblock The spaces of index one minimal surfaces and constant mean curvature
  surfaces embedded in flat three manifolds.
\newblock {\em Trans. Amer. Math. Soc.}, 348:391--410, 1996.

\bibitem{RosMSRI}
Antonio Ros.
\newblock The isoperimetric problem.
\newblock Lecture notes on a series of talks given at the Clay Mathematics
  Institute Summer School on the Global Theory of Minimal Surfaces June 25,
  2001 to July 27, 2001, at the Mathematical Sciences Research Institute,
  Berkeley, California, {\tt http://www.ugr.es/\~{}aros/isoper.pdf}.

\bibitem{Schwarz}
H.~A. Schwarz.
\newblock {B}eweis des {S}atzes, dass die {K}ugel kleinere {O}berfl\"ache
  besitzt als jeder andere {K}{\"o}rper gleichen {V}olumens.
\newblock {\em Nachrichten K{\"o}niglichen Gesellschaft Wissenschaften
  G{\"o}ttingen}, pages 1--13, 1884.

\bibitem{Simon}
Leon Simon.
\newblock {\em Lectures on Geometric Measure Theory}.
\newblock Proc. Centre Math. Anal, Australian Nat. U. Vol. 3. Centre for
  Mathematical Analysis, Australian National University, Australia, 1984.

\bibitem{Taylor}
Jean Taylor.
\newblock The structure of singularities in soap-bubble-like and soap-film-like
  minimal surfaces.
\newblock {\em Ann. Math}, pages 489--539, 1976.

\bibitem{WP}
D.~Weaire and R.~Phelan.
\newblock A counter-example to {Kelvin's} conjecture on minimal surfaces.
\newblock {\em Phil. Mag. Lett.}, 69:107--110, 1994.

\bibitem{foams}
Denis Weaire and Stefan Hutzler.
\newblock {\em The Physics of Foams}.
\newblock Oxford University Press, Oxford, 2001.

\end{thebibliography}
\bibliographystyle{plain}
\end{document}